\documentclass[11pt]{amsart}
\usepackage{geometry}                
\usepackage{graphicx}
\usepackage{amssymb}
\usepackage{epstopdf}
\title{Undefinability and Absolute Udnefinability in Arithmetic}
\author{Roman Kossak}
\newtheorem{thm}{Theorem}[section]
\newtheorem{cor}[thm]{Corollary}

\newtheorem{prop}[thm]{Proposition}
\newtheorem{defn}[thm]{Definition}

\def\vp{\varphi}
\def\al{\alpha}
\def\om{\omega}

\def\fa{\forall}
\def\ex{\exists}
\def\iff{\Longleftrightarrow}
\def\then{\Longrightarrow}
\def\into{\longrightarrow} 

\def\elem{\el\  extension}
\def\ar{arithmetic}
\def\ct{countable}
\def\au{automorphism}
\def\aug{automorphism group}

\def\cont{$2^{\aleph_0}$}

\def\st{structure}
\def\fo{first-order}

\def\nn{natural number}

\def\be{\begin{xca}}
\def\ee{\end{xca}}

\def\elem{elementary extension}
\def\iso{isomorphic}
\def\isom{isomorphism}

\def\rs{recursively saturated}

\def\lom{L_{\om_1,\om}}
\def\pt{pointwise}
\def\pd{parametrically definable}
\def\lpa{\lan_{\sf ar}}
\def\resp{resplendent}

\def\N{{\mathbb{N}}}
\def\Z{{\mathbb{Z}}}

\def\A{{\mathfrak A}}

\def\M{{\mathfrak M}}

\def\lan{{\mathcal L}}

\def\PA{{\sf PA}}
\def\Pr{{\sf Pr}}

\def\zfc{{\sf ZFC}}

\def\TO{$\Th(\N,<)$}
\def\TS{$\Th(\N,S)$}

\def\Tr{{\sf Tar}}

\DeclareMathOperator{\Th}{Th}

  \newcommand{\god}[1]{\left\ulcorner#1\right\urcorner}

\begin{document}
\maketitle

\begin{abstract} This is a survey of results on definability and undefinability in models of \ar. The goal is to present a stark difference between undefinability results in the standard model and  much stronger versions about expansions of nonstandard models. The key role is played by counting the number of automorphic images of subsets of \ct\ \resp\ models of Peano Arithmetic.
\end{abstract}

\section{Introduction}

By a  {\it language} of a structure we mean its set of  function, relation, and constant symbols.  The {\it language of \ar} is $\{+,\times\}$.

Let $\N$ be the set of \nn s. All computable and computably enumerable sets of \nn s are \fo\ definable\footnote{In this note definability will mean definability without parameters. If parameters are involved we will refer to parametric definability. For the standard model both notions coincide.} in the standard model $(\N,+,\times)$.  Beyond that there is an infinite hierarchy of definable  sets whose  complexity is measured by the number of alternations of quantifiers in their definitions.   First-order logic is strong enough to capture all this complexity.   This was first revealed  in  the proof of G\"odel's   incompleteness theorems. G\"odel showed how recursive definitions can be converted to \fo\ ones, which opened the door to an  interpretation of the syntax of \fo\ logic in  \ar\ (arithmetization). Under this interpretation, each \fo\ \ar\ formula $\vp$  is assigned a code,  its G\"odel number, $\god{\vp}$. It can be shown that  for each \nn\ $n$, the set  ${\rm Tr}_n$  of G\"odel numbers of sentences with at most $n$ quantifiers which are true in the standard model is definable.  The defining formulas follow Tarski's definition of truth, but still it takes an effort to write them down explicitly---especially in the case of $n=0$,  and to show in Peano Arithmetic (\PA) that they have the required property. All the gory details are given in \cite[Chapter 9]{kaye}.

Soon after G\"odel results,  Skolem's construction of an elementary extension of the standard model exposed a weakness of \fo\ logic by showing that the standard model of arithmetic is not uniquely determined by its \fo\ theory. Almost at the same time,  Tarski proved a general result on formal theories from which it follows that the ``full truth"  ${\rm Tr}=\bigcup_{n\in \N}{\rm Tr}_n$ is not definable in the standard model. Both results were seminal for later  developments in mathematical logic.

 The question of definability becomes  important in the case of functions and relations that may be regarded as intrinsic to the \st. In model theory, if a \st\ is {\it given}, then  what is given  are its functions, relations, and constants. It does not matter how they are given; they may be  explicitly defined, but still highly complex; they  may be  given only by  existence theorems of set theory, or they may  just be assumed to exist. These sets and constants  define the \st\ explicitly. What is also implicitly given is the set of all \fo\ definable subsets of all finite Cartesian powers of the domain.  To know a \st\ is to know its definable sets. First-order definitions reveal the complexity of the sets they define. They tell us how the definable sets are built from the basic relations, functions, and constants by means of Boolean operations, Cartesian products, and projections.  First-order logic gives us what I call {\it logic visibility}. We {\it see} the geometry of the definable sets  through the eyes of logic. 

The fact that  truth is undefinable is a deficiency of \fo\ logic. The set Tr is the union of countably many definable sets. Moreover, while the quantifier complexity of the sets ${\rm Tr}_n$ increases with $n$, the definitions form  a well-defined recursive sequence. A number is a G\"odel code of a true sentence if and only if  it is in one of the sets defined by these formulas, and this is an example of a perfectly clear mathematical definition. In other words, Tr is an example of an intrinsic set which \fo\ logic does not see.  

Truth becomes  logically visible in an extension of \fo\ logic known as $\lom$.   In $\lom$, in addition to  all \fo\ operations, we can form conjunctions and disjunctions of arbitrary \ct\ sequences of formulas, as long as the free variables of all formulas in the sequence are contained in a finite set. This extends definability to \ct\ intersections and \ct\ unions of definable subsets of a fixed Cartesian power of the domain of a \st. The extension is radical. The following theorem was  proved independently by David Kueker \cite{kueker} and Gonzalo Reyes \cite{reyes}.

\begin{thm}\label{kuerey} Let $\A$  be a   \ct\ \st\  for a \ct\ language. Then,  a relation $R$ on $\A$ has an $\lom$ definition with finitely many parameters   if and only if $R$ has at most countably many  images under \au s of $\A$. Moreover, if the set of automorphic images of $R$ is uncountable,  it must be of power continuum.
\end{thm}

 The Kueker-Reyes theorem is a strengthening of an earlier  result  of Dana Scott \cite{scott}. Scott proved that every relation on the domain of a \ct\ \st\  is fixed setwise by every \au\ if and only if has an $\lom$ definition without parameters. In particular, in every rigid \ct\ structure  every element of the domain has an $\lom$ definition, and it follows form the Kueker-Reyes theorem that in every \ct\  \st s with fewer than \cont\ \au s, every relation has an $\lom$ definition, possibly with parameters. 
 
 One direction of Kueker, Reyes, and Scott's results is obvious. Parametrically definable sets are setwise fixed by \au\ that fix all parameters in their definitions. This applies not only to \fo\ logic and $\lom$, but in fact it is a requisite for all logical formalisms. Logical properties should be preserved by \isom s. Thus, if a relation in a \ct\  \st\ has uncountably many  automorphic images, I will call it {\it absolutely undefinable}. 

  Let $S$ be the successor relation on $\N$.  In the sequence $(\N,S)$, $(\N,<)$, $(\N,+)$, $(\N,+,\times)$, $(\N,+,\times, {\rm Tr})$,  each \st\ is richer than the previous one, for it has a new relation or function that is not definable in the previous one. This is telling us something about logic, but not that much about the nature of those functions and relations, as they all are $\lom$-definable in $(\N,S)$ by simple formulas that correspond to the first-order recursive definitions. In contrast,  absolute undefinability theorems about expansions of nonstandard models, which we are going to discuss,  reveal strong independence of  the  basic arithmetic functions and relations in nonstandard models.  
First, we will examine in detail definability of linear orderings in \elem s of $(\N,S)$.  This will be  followed by a number of absolute undefinability results about expansions of \ct\ models of Presburger \ar\ and expansions of \ct\ \resp\ models of \PA\ to axiomatic fragments  of the \fo\ theory of $(\N,+,\times, {\rm Tr})$. 

\thanks{ Ali Enayat,  Mateusz {\L}e{\l}yk, and Simon Heller  have  reviewed  the draft of this paper and provided valuable comments. I also want  to thank Alfred Dolich, Emil Je\v{r}\'{a}bek, Zachiri McKenzie, and Simon Heller for their observations  included in sections 3 and 4.}

\subsection{Resplendence}
For a \st\ $\A$,   $\Th(\A)$ is  the complete \fo\ theory of $\A$. A \st\ $\A$ is  {\it resplendent} if for all for all tuples $\bar a$ in the domain of $\A$ and all sentences $\vp(R,\bar a)$ in the language of $\A$ with an extra relation symbol $R$, if  some model  of $\Th(\A,\bar a)$   is expandable to a model of $\vp(R,\bar a)$, then already  $\A$ is expandable to model of $\vp(R,\bar a)$.  

Definition of resplendence does not require any assumptions about the language, but the results that we will discuss do. {\bf From now on the will assume that the language (signature) of each \st\ is finite.} 

By a theorem proved by Barwise and Schlipf, and independently by Ressayre, a \ct\ \st\ is resplendent if and only if it is \rs. Every \ct\ \st\ has a \ct\ \resp\  \elem. Moreover, \ct\ \resp\ models can be characterized by the following stronger property:  Let $\A$ be a \ct\ resplendent \st,   let $\lan$ be the language of $\A$, and let $T$ be a computably axiomatized theory in a  language $\lan'\supseteq \lan\cup\{\bar a\}$, for  $\bar a$ in the domain of $\A$. If some model of $\Th(\A,\bar a)$  is expandable to a model  of $T$, then $\A$ has an expansion to a \resp\ model of $T$. This property  is called \emph{chronic resplendence}.

We will say that an expansion of a \st\ $\A$ is {\it (parametrically) definable} if all new functions, relations, and constants are (parametrically) definable in $\A$. If an \elem\ of a \st\ $\A$ has  a (parametrically) definable expansion to a model of $\vp(R,\bar a)$, then  $\A$ has such an expansion.  
 
For a brief introduction to resplendent models see \cite{resplendent}, and for a full discussion of the role of resplendence in models of \ar\ see \cite{smo} and \cite[Chapter 15]{kaye}. For us, the following result of Schlipf \cite{schlipf} is crucial. 

\begin{thm}\label{crucial}  If $(\A,R)$ is   \ct\ and \resp,  and   $R$  is not \pd\ definable in $\A$,  then $R$ is absolutely undefinable in $\A$.
\end{thm} 

Suppose that  $T$ is  a computably axiomatizable theory in the language of a \ct\ \resp\ \st\ $\A$ with an additional relation symbol,  $T$ is  consistent with $\Th(\A)$ and  $\A$ has has an expansion to a model of $T$ that is not \pd.  Then, by chronic resplendence, $\A$ has a  \resp\ expansion $(\A,R)$ to a model of $T$ that is not \pd. It follows from Schlipf's theorem that each such $R$ is absolutely undefinable in $\A$. This shows that absolute undefinability is hard to avoid. If a \ct\ resplendent model has expansions to a model of $T$ that are not \pd, then among those expansions there have to be absolutely undefinable ones.  

A classification of \ct\ models of \TS\ is given in the next section. All statements about the following example will be justified there.  

Let  $I$ be a unary relation symbol, and let $\vp(I)$ be the sentence 
\[[\fa y\lnot S(y,x)\then I(x)]\land [\fa x,y (I(x)\land S(x,y))\then I(y)]\land \ex x\lnot I(x).\]
Clearly, $(\N,S)$ is not expandable to a model of $\vp(I)$, but  every \elem\ of $(\N,S)$ is.  This shows  that $(\N,S)$ is not \resp. Because $(\N,S)$ satisfies the induction schema, no expansion of a model of \TS\ to a model of $\vp(I)$ can be \pd, but  we will see that some such expansions are $\lom$-definable with parameters;  hence they are  not absolutely undefinable. By Schlipf's theorem, $(\N,S)$ has an \elem\ to a \ct\ \resp\ model which has an absolutely undefinable expansion to a model of $\vp(I)$.  We will see that \TS\ has exactly one such model.

\subsection{Real and Imaginary sets}
The last section of this paper is devoted to some specific results about intrinsic but absolutely undefinable sets in \ct\ \resp\ models of \PA. Similar results were obtained independently by Athanassios Tzouvaras \cite{tzou}. Tzouvaras uses the real/imaginary terminology introduced by \v{C}uda and Vop\v{e}nka in the context of Alternative Set Theory \cite{cuvo}. If $\M$ is a countable \resp\ model  a relation  $R$ on   $\M$ is real if and only if the set of automorphic images of $X$ is countable.   Thus, all imaginary relations are absolutely undefinable. Among many interesting results, Tzouvaras showed that all nontrivial \au s of \ct\ \resp\ models of \PA\ are imaginary. This was also shown independently by another proof in \cite{kkk}. Both proofs are short, but they use nontrivial results about models of \PA: Tzouvaras uses Ehrenfeucht's Lemma, the authors of \cite{kkk} use Kotlarski's Moving Gaps Lemma.

\section{The successor relation}

There is hardly a simpler \st\ than $(\N,S)$, but despite its simplicity, the addition and multiplication are uniquely determined in it and  so is the set Tr. For all numbers $m$ and $n$, $m+0=m$ and $m+(n+1)=(m+n)+1$. It is a definition by recursion. To compute $m+n$ we can  start with  $m$ and move $n$ successors up. Why is it that we see  clearly something that the powerful \fo\ logic cannot?  Our understanding of $(\N,S)$ includes the fact that every \nn\ can be reached from 0 by a {\it finite} number of successor steps. It is this finiteness that \fo\ cannot handle.  For each $n$, the formula  
\[\vp_n(x,y)=\ex x_1\ex x_2 \cdots \ex x_n[x_1=x\land S(x_1,x_2)\land \cdots \land x_n=y]\]
defines  the relation $x+n=y$ in $(\N,S)$. However,  there is no \fo\ formula that defines $x+z=y$. There is no \fo\ way to say ``repeat \dots\ $z$-times."  How do we know this? A  proof will be given below.  

A \st\ is {\it\pt\ definable} if each element of its domain has a \fo\ definition.
 In $(\N,S)$, the formula $\sigma_0(x)=\fa y \lnot S(y,x)$ defines 0, and for each $n>0$, $\sigma_n(x)=\vp_n(0,x)$ defines $n$; hence $(\N,S)$ is \pt\ definable, and it follows that every model of \TS\ has an elementary submodel which is isomorphic to $(\N,S)$. We will identify this submodel with $(\N,S)$ and call its elements {\it standard}; all other elements will be called {\it nonstandard}. A model is called nonstandard if it has nonstandard elements. The same terminology will be applied to models of  theories of expansions of $(\N,S)$. {\bf The formulas $\vp_n(x,y)$ and $\sigma_n(x)$ are fixed for the rest of the article.}

Because in $(\N,S)$ every element other than 0 has a successor and a predecessor, it follows that  every  nonstandard element $c$ in an \elem\ of $(\N,S)$ belongs to a chain which is isomorphic to the successor relation on the set of integers. We will denote such a chain by $[c]$ and call it the $\Z$-chain of $c$.

By the compactness theorem, there is a model of \TS\ with at least two $\Z$-chains. Let $(M,S)$ be such a model. We will use $(M,S)$ to show that no linear ordering is \pd\ in $(\N,S)$. Because $(\N,S)$ is \pt\ definable, it is enough to consider formulas without parameters. To get a contradiction, suppose $\vp(x,y)$ defines  a linear ordering of $\N$. Then $\vp(x,y)$  also defines a linear ordering $\prec$ in $(M,S)$. Let $[c]$ and $[d]$ be distinct $\Z$-chains in $M$ with $c\prec  d$, and let $f(c+n)=d+n$ and $f(d+n)=c+n$, for all $n\in \Z$,  and $f(x)=x$ for all other $x$ in $M$. Then $f$ is an \au\ of $(M,S)$, and $f(d)\prec f(c)$, which shows  that $\prec$ cannot be definable in $(M,S)$ giving us a contradiction. In particular, the usual ordering $<$ is not definable in $(\N,S)$.  Because  $S$ is definable in $(\N,<)$, it follows that $(\N,<)$ is a proper expansion of $(\N,S)$. 

In every \pt\ definable model every relation has a parameter-free $\lom$ definition. It is just the elementary diagram of the relation coded by a single $\lom$ sentence.  For example, the ordering of  $\N$ is defined by $\bigvee\{\sigma_m(x)\land\sigma_n(y): m<n\}$. In this case there is also a  definition which ties the ordering to the successor relation in a straightforward way: $x<y$ if and only $\bigvee_{n>0} \vp_n(x,y)$.

Using Ehrenfeucht-Fra\"iss\'e game characterization of first-order elementary equivalence, it can be shown that the extension of $(\N,S)$ by any number of  $\Z$-chains is  a model of \TS\  (see \cite[Section 2.4]{marker}). Any two models with the same number of $\Z$-chains are isomorphic, so up to isomorphism there are exactly $\aleph_0$ \ct\ models of \TS. Except for the standard model, none of those models is rigid. Each \au\ of a nonstandard model of \TS\ is a composition of a permutation of the $\Z$-chains and an \au\ which  for each  chain either fixes it \pt\ or shifts its elements either up or down. It follows that  \aug s of  models with finitely
many $\Z$-chains are \ct, and the \aug\ of the model with $\aleph_0$ chains is of power \cont. Each  model of \TS\ can be expanded by a linear ordering that is compatible with the successor relation.We can first order the $\Z$-chains,  and then order the elements in each chain according to the successor relation. In the case of models with finitely many chains, different orderings of the chains, give us different, but isomorphic  expansions; hence, each ordering has finitely many automorphic images, and---by the Kueker-Reyes theorem---they must have $\lom$ definitions with parameters. Such definitions can be obtained by choosing one parameter from each $\Z$-chain, deciding how they are ordered, and then defining the ordering of each $\Z$-chain by a formula  with the chosen parameters similar to the one we gave  for the ordering of the standard model.

In the case of models with countably many $\Z$-chains, it is not hard to see  that each  ordering  compatible with the successor relation has \cont\ automorphic images; hence it is absolutely undefinable. 

Resplendent models of \TS\  have infinitely many $\Z$-chains. For example, if $(M,S)$ has only one $\Z$-chain, let $a$ be any element of the chain and consider the following theory $T$ with a unary relation symbol $A$: 
\[\{\ex ! x A(x)\land [\fa x [A(x)\then[\lnot \sigma_n(x)\land \lnot \vp_n(a,x)\land \lnot \vp_n(x,a)]: n\in \N\}.\] 
By the compactness theorem, $(M,S,a)$ has an elementary extension that is expandable to a model of $T$. Because $(M,S)$ has no such expansion; hence it is not \resp.   This argument can be generalized to models with finite numbers of $\Z$-chains.  Hence, up to \isom,   there is only one \ct\ \resp\ model of \TS.

\section{The ordering and the additive structure}

 The discussion from the previous section shows that the nonstandard part of each model \TO\ is the union of its  $\Z$-chains.  The set of $\Z$-chains of a model of \TO\ is linearly ordered by the ordering inherited from the model. Using Ehrenfeucht-Fra\"iss\'e games, it can be shown  that for any linearly ordered set $(I,<)$ the union of $(\N,<)$  and the ordered set of $\Z$-chains which is isomorphic to $(I,<)$ is a model of \TO. Because there are \cont\ nonisomorphic \ct\ linearly ordered sets, it follows that \TO\ has  \cont\ nonisomorphic \ct\ models. 

We will now show that addition is not definable in $(\N,<)$.  We will prove a stronger fact: $(\N,<)$ is {\it minimal}, i.e., every  subset of $\N$ definable in $(\N,<)$ is either finite or cofinite.\footnote{In general, a \st\ is minimal if every \pd\ subset of its domain is either finite or cofinite. Because $(\N,<)$ is \pt\ definable, we only consider formulas without parameters.} To get a contradiction suppose that the set defined in $(\N,<)$ by $\vp(x)$ is neither finite nor cofinite. Then\footnote{We can use $S(y,z)$ in the formula below, because $S$ is definable in $(\N,<)$.}
\[(\N,<)\models \fa x \ex y \ex z[x<y\land S(y,z)\land \vp(y)\land \lnot\vp(z)].\eqno{(*)}\]
Let $(M,<)$ be an \elem\ of $(\N,<)$.   Then by $(*)$ there must be  nonstandard $c$  and $d$ such that 
\[(M,<)\models S(c,d)\land \vp(c)\land\lnot\vp(d).\]
 Let $f(x)$ be the successor of $x$ for all $x$ in the $\Z$-chain $[c]$,  and let $f(x)=x$ for all other $x$.  Then $f$ is an \au\ of $(M,<)$, and $f(c)=d$, which contradicts the fact that $(M,<)\models \vp(c)\land \lnot\vp(d)$ and finishes the proof.

$(\N,+)$ is not minimal. For example, the  formula $\ex y [x=y+y]$ defines the set of even numbers. It follows that addition is not definable in $(\N,<)$. Because the usual ordering of the \nn s is definable in $(\N,+)$, it follows that $(\N,+)$ is a proper expansion of $(\N,<)$.

We will denote $\Th(\N,+)$  by \Pr\ (Presburger Arithmetic). In the previous section we saw that every model  of \TS\ is expandable to a model of \TO. Now we will briefly discuss expandability of models of \TO\ to models of \Pr. 

For every natural number $n$, either $n$ or $n+1$ is even, and this can be expressed by a \fo\ sentence. Hence, every $\Z$-chain of a model of \Pr\ can be presented as $[a]$ for an even $a$.  Let $a$ be  nonstandard and even.  It is easy to prove that if $a=b+b$, then $b$  is nonstandard and  $[b]$, $[a]$, and $[a+a]$ are disjoint. It follows that in a nonstandard model of \Pr\ there is no smallest and no largest $\Z$-chain. If $a$ and $b$ are even, $a<b$, $[a]$ and $[b]$ are disjoint, and $c+c=a+b$, then $a<c<b$ and $[a]$, $[c]$, and $[b]$ are disjoint. This shows that the ordered set of $\Z$-chains in any \ct\ nonstandard model of \Pr\ is isomorphic to the ordered set of rational numbers; hence, up to isomorphism, there is only one \ct\ model \TO\ which is expandable to a model of \Pr. At this point, model theory becomes harder. It can be shown that there are \cont\ isomorphism types of such expansions, but we will not discuss the details here. 

Again resplendence enters the picture. Up to isomorphism, there is only one \ct\ \resp\ model of \TO. It is the model whose ordered set of $\Z$-chain is \iso\ to the ordered set of rational numbers, so a \ct\ model of \TO\ is expandable to a model \Pr\ if an only if it is \resp.  The question now is: How many such expansions are there and are there any that are not absolutely undefinable? Because there are  \cont\ \isom\ types of \ct\ models of \Pr, it follows that the \ct\ \resp\ model of \TO\ can be expanded to a model of \Pr\ in \cont\ nonisomorphic ways.

To expand a model $(M,<)$ of \TO\ to a model of \Pr\ whose addition is compatible with the ordering, we need a function $f:M^2\into M$  such that $(M,f)\models \Pr$ and \[(M,<,f)\models \fa x,y[(\lnot(x=y)\land \ex z f(x,z)=y))\iff\ x<y].\]
The required property of the expansion is no longer a single sentence, but an infinite theory. By a theorem of Presburger, \Pr\ is computably axiomatized; hence resplendence can still be applied. Let $\M$ be the \ct\ \resp\ model of \TO. It follows from  Schlipf's theorem that $\M$ has an absolutely undefinable expansion to a model of \Pr.  In a response to my query Emil Je\v{r}\'{a}bek gave a short argument proving that in fact each expansion of $\M$ to a model of \Pr\ is absolutely undefinable (unpublished).

While there are no  nonstandard rigid models of \TO, there are \ct\ nonstandard rigid models of \Pr, but they are rare.  For a full discussion see \cite{jerabek}.  

\subsection{Membership and inclusion}
The relationship between the addition and the ordering of models of \Pr\ is similar to that of the relations membership and inclusion of models of set theories. In \cite{hamkik} and a follow-up paper \cite{hamkik2}, Joel David Hamkins and Makoto Kikuchi showed that the membership relation is not definable in inclusion reducts $(M,\subseteq^M)$ of models $(M,\in^M)$ of \zfc\ and other set theories. Moreover, if $(M,\in^M)$ is a \ct\ model of \zfc, then $(M,\subseteq^M)$ is $\om$-saturated---hence it is \resp---and up to isomorphism there is only one such reduct. Zachiri McKenzie  noticed that using the results of Hamkins and Kikuchi one can show that in each inclusion reduct of a \ct\ model of \zfc, $\in^M$ is absolutely undefinable.

\section{The multiplicative structure}

A theorem of Ginsburg and Spanier  says that the subsets of $\N$ which are definable in $(\N,+)$ are exactly the ultimately  periodic ones, i.e., the sets $X$ for which there is a $p$ such that for sufficiently large $x$, $x\in X$ if and only if $x+p\in X$ \cite{ginspa}. The set of squares is not ultimately periodic, but it is definable in $(\N,\times)$; hence multiplication is not definable in $(\N,+)$. In the other direction, addition is not definable in $(\N,\times)$. There are many arguments which can be used to show this. A short one uses the fact that $(\N,\times)$ admits lots and lots of \au s. Every permutation $\al$ of the set prime numbers can be extended to an \au\ $f$ of $(\N,\times)$ defined for each product of prime numbers by $f(\Pi_{i<n}p_i^{k_i})=\Pi_{i<n}\al(p_i)^{k_i}$.  This shows that even the successor relation is not definable in $(\N,\times)$. Incidentally, this is the only obstacle.  Addition is definable in $(\N,\times, S)$. It is easy to check that the following holds for all  $x$, $y$, and all $z\not=0$.
\[x+y=z \iff (zx+1)(zy+1)=z^2(xy+1)+1.\]

Model theory of the natural numbers becomes more difficult when addition and multiplication are joined together.  What made the previous discussion easier was the fact that the \fo\ theories of $(\N,S)$, $(\N,<)$, and $(\N,+)$ are computably axiomatized, and so is $\Th(\N,\times)$, but this is more difficult to show. An axiomatizations  for $\Th(\N,\times)$ was first given  by Patrick C\'egielski in 1981.\footnote{For references and a comprehensive survey of axiomatizability of \fo\ theories of number theoretic structures see \cite{salehi}.}  

The first-order theory of $(\N,+,\times)$ is not computably axiomatizable.  To be able to apply model theory of resplendent \st s,  we need to shift to  axiomatizable subtheories. We will focus on \PA. Most of the results we will mention apply to fragments of \PA\ as well as to its extensions, most notably to $\PA^*$, which is Peano axioms in any \ct\ language extending the language of \ar.

An important result in model theory of \ar\ is that if $(M,+,\times)$ is a nonstandard model of \PA, then $(\M,+)$ is a model of \Pr, $(M,\times)$ is a model of $\Th(\N,\times)$,  and both  $(M,+)$ and $(M,\times)$ are \rs; hence, if $M$ is \ct, then both resducts are resplendent. Moreover, if $(M,+,\times)$ and $(N,+,\times)$ are \ct\ nonstandard models of \PA, then $(M,+)$ and $(N,+)$ are \iso\ if and only if $(M,\times)$ and $(N,\times)$ are.   For a discussion of these results and a proof that the equivalence above no longer holds in the un\ct\ case see \cite{kns}.

Every \ct\ \resp\ model  of \Pr\ can be expanded to \cont\ non\iso\ models of \PA, including \cont\ rigid models, \cont\ resplendent models, and \cont\ other models, and each of these expansions is absolutely undefinable. All these statements, except the last one follow from the standard results about \ct\ \resp\ models of \Pr\ summarized in \cite[Chapter 15]{kaye}. The last  result is due to Alf Dolich and Simon Heller (unpublished). It is based on the fact that every \resp\ model of \Pr\ has \au s which move some elements to their squares. In the special case when   $(M,+, \times)$ is a \ct\ \resp\ model of \PA,  it follows directly from Schlipf's theorem that $\times$ is absolutely undefinable in $(M,+)$, and that $+$ is absolutely undefinable in $(M,\times)$. 

\section{Tarski Arithmetic}
We begin this section with a general undefinability result, from which Tarski's theorem on undefinability of truth follows.

\begin{defn} Let $M$ be the domain of a model $\M$. A set $U\subseteq M^2$ is \emph{universal} if for each \pd\ $X\subseteq M$  there is an element  $b\in M$ such that
$X=\{a: (a,b)\in U\}$.
\end{defn}

\begin{prop} No universal set is \pd.
\end{prop}

\begin{proof} Suppose that a \pd\ $U\subseteq M^2$ is universal. Let $D=\{a: (a,a)\notin U\}$. $D$ is \pd\, so there is $b$ such that $D=\{a: (a,b)\in U\}$.
We have:  $b\in D$ iff $(b,b)\notin U$ iff $b\notin D$. Contradiction.
\end{proof}
In the corollary below we will take advantage of \ar\ coding of finite sequences, which allows us to reduce definability with parameters to definability with only one parameter. In  the formula below, $\langle x,y\rangle$ denotes Cantor's pairing function, and $\god{\vp(x)}$  is the G\"odel number of $\vp(x)$.
\begin{cor}[Undefinability of Truth]\label{tarski}  Let $\M$ be a model of \PA.  There are no \pd\ $S\subseteq M^2$  such that for all \fo\ formulas $\vp(x)$ of the language of \ar\ with parameters from $M$ and all $a\in M$ \[(a, \god{\vp(x)})\in S \iff  \vp(a).\] 
\end{cor}
In the corollary, it does not matter how  $\god{\vp(x)}$ is defined, all we need is that  it is an element of $M$.

It follows from Corollary \ref{tarski} that $(\N,+,\times, {\rm Tr})$ is a proper expansion of $(\N,+,\times)$.  Let \Tr\ (Tarski Arithmetic) be $\Th(\N,+,\times, {\rm Tr})$.  A routine argument using the overspill principle shows that if $(\M,T)$ is a nonstandard model of \Tr, then $\M$ is \rs; hence, if $\M$ is \ct, then it is resplendent. 

By G\"odel's incompleteness theorem, \Tr\ is not axiomatizable. In recent years, much attention was given to the study of axiomatic subtheories of \Tr. See \cite{cies} for a comprehensive account. Originally, the results were formulated in terms of satisfaction classes (binary), as we will do below. Recently, the emphasis has shifted to truth predicates (unary). With a modicum of induction in the axiom system, satisfaction classes and truth predicates are definable from one another. In the absence of induction, the situation is more subtle, see \cite{cies2} for details.  Let Sat be a binary relation symbol. We will consider the following theories in the language $\lpa({\rm Sat})=\{+,\times, {\rm Sat}\}$. 
\begin{enumerate}
\item {\sf Partial satisfaction class:} \PA\ + all sentences of the form \[\fa y[\vp(y)\iff {\rm Sat}(\god{\vp(x)},y)].\] 
\item {\sf Partial inductive satisfaction class:} (1) + induction for all formulas of $\lpa({\rm Sat})$.
\item{\sf  Full satisfaction class:} Arithmetized axioms for satisfaction (but no induction). For example, the axiom of fullness: \[\fa \vp \fa z [{\rm Form}_\PA(\vp)\then ({\rm Sat}(\lnot\vp,z))\iff \lnot {\rm Sat}(\vp,z))],\]
where ${\rm Form}_\PA(x)$ is an \ar\ formula representing in \PA\ the set of G\"odel numbers of \ar\ formulas with one free variable. 
\item {\sf Full inductive satisfaction class:} (3) + induction for all formulas of $\lpa({\rm Sat})$.
\end{enumerate}
By Corollary \ref{tarski}, no partial satisfaction class is \pd. 

  It was shown in \cite{koskot1} that all partial inductive partial satisfaction classes in \ct\ \resp\ models of \PA\ are absolutely undefinable. This turned out to be a special case of a much more general result, for which we need a couple of definitions. A subset $X$ of  a model of \PA\ is a {\it class} if for every $a$, $\{x\in X: x<a\}$ is parametrically definable. It is not difficult to show that $X\subseteq M$ is a class of a model $\M$ of \PA\ if and only if $(\M,X)$ is a model of bounded induction $I\Delta_0(X)$.  If $(\M,X)$ is a model of $\PA^*$, we call $X$ {\it inductive}.  All inductive sets are classes. Every \ct\ model of \PA\ has classes which are not inductive, and also has inductive sets which are not \pd. There are un\ct\ models of \PA\ all of whose classes are definable (the rather classless models). 
 
 Let $\M$ be a  \ct\ \resp\ models of \PA. For such models,  Jim Schmerl proved that  all undefinable classes  of $\M$ are absolutely undefinable \cite[Theorem 8.2.3]{tsomopa}. This shows that if $S$ is a partial inductive satisfaction class and $(\M,S)\models I\Delta_0(S)$, then $S$ is absolutely undefinable.   In the absence of any amount of induction, the case of full satisfaction classes requires a separate treatment. By a theorem of Kotlarski, Krajewski, and Lachlan, $\M$ has a full satisfaction class,   but by a theorem of Bartosz Wcis{\l}o \cite{wcilel}, $\M$ can have  a full satisfaction class that is a class only under certain assumptions about $\Th(\M)$. Combining several results one can show that all full satisfaction classes in \ct\ \resp\ models are absolutely undefinable. A direct proof of this result and a comprehensive discussion are in  the recent paper \cite{bart}.

\bibliographystyle{amsalpha}
\bibliography{absolute}

\end{document}